\documentclass[a4paper,pdflatex,sn-mathphys-num]{sn-jnl}

\usepackage{graphicx}
\usepackage{multirow}
\usepackage{amsmath,amssymb,amsfonts}
\usepackage{amsthm}
\usepackage[title]{appendix}
\usepackage{xcolor}
\usepackage{booktabs}
\usepackage{url}
\usepackage{algorithm}
\usepackage{float}
\usepackage{algpseudocode}

\theoremstyle{thmstyleone}

\theoremstyle{thmstyletwo}

\theoremstyle{thmstylethree}

\begin{document}

\title[Sigmoid Tucker Decomposition]{ Sigmoid Tensor Decomposition for Compact Representations of Binary Data}

\author*[1, 2]{\fnm{Irina} \sur{Kopnina}}\email{kopninaia@my.msu.ru}
\author[1, 2]{\fnm{Sergey} \sur{Matveev}}
\author[1, 3]{\fnm{Ivan} \sur{Oseledets}}


\affil[1]{\orgname{Lomonosov Moscow State University(MSU)}, \orgaddress{\city{Moscow}, \country{Russia}}}
\affil[2]{\orgname{Marchuk Institute of Numerical Mathematics, Russian Academy of Sciences (INM RAS)}, \orgaddress{\city{Moscow}, \country{Russia}}}
\affil[3]{\orgname{Artificial Intelligence Research Institute (AIRI)}, \orgaddress{\city{Moscow}, \country{Russia}}}

\abstract{Nonlinear low-rank decompositions are receiving growing attention as a way to represent structured data with fewer parameters than conventional linear methods. Existing nonlinear approaches have focused mainly on matrices, which limits their ability to capture multidimensional structure. We propose a nonlinear tensor decomposition based on sigmoid nonlinearization in the Tucker and tensor-train formats. A binary tensor is represented by a low-rank real-valued score tensor, whose parameters are estimated by minimizing the binary cross-entropy loss. The resulting optimization problem is jointly nonconvex, but becomes convex with respect to one factor or score block when the remaining factors are fixed, as in the corresponding matrix case. We evaluate the method on several structured binary tensors and compare it with conventional Tucker and tensor-train decompositions across a range of common ranks. These results indicate that sigmoid nonlinearization is particularly effective in the low-rank regime, where it can provide more accurate binary reconstructions than conventional linear tensor decompositions at comparable model sizes.
}



\maketitle

\section{Introduction}\label{sec:intro}

Low-rank approximation is one of the essential tools for representing large data arrays with a small number of parameters.
Such representations are widely used in scientific computing, machine learning, and other data-intensive applications. Many multidimensional datasets exhibit correlations across modes that allow them to be represented by low-dimensional models, thereby reducing storage requirements and computational costs.

For a matrix $\mathbf{A}\in\mathbb{R}^{M\times N}$ and a target rank $r\ll\min(M,N)$, classical low-rank approximation seeks factors $\mathbf{U}\in\mathbb{R}^{M\times r}$ and $\mathbf{V}\in\mathbb{R}^{N\times r}$ by solving
\[
\min_{\mathbf{U},\mathbf{V}}\|\mathbf{A}-\mathbf{U}\mathbf{V}^{\mathsf T}\|_F^2.
\]
The Frobenius norm is particularly important because the Eckart--Young--Mirsky theorem provides an explicit optimal rank-$r$ approximation in this norm~\cite{eckart1936,mirsky1960}. Related optimality results also hold for  other unitarily invariant norms, including the spectral norm~\cite{von1937some, mirsky1960}. 

For discrepancies such as entrywise $L_1$ loss or Kullback--Leibler divergence, no comparable universal singular-value formula is generally available. Nevertheless, numerical methods can efficiently compute useful low-rank representations under these objectives~\cite{chierichetti2017l1,lee2001nmf}. In addition to selecting a discrepancy, one may impose constraints on the factors themselves when the application provides structural or physical prior information. Nonnegative matrix factorization, sparse factorization, and related constrained models encode such information through the feasible set or through explicit regularization terms~\cite{lee1999nmf,hoyer2004sparse,mairal2010online}. 

More generally, the discrepancy should be compatible with the assumed observation model: the squared Frobenius loss corresponds to the Gaussian least squares, whereas for binary observations the Bernoulli negative log-likelihood gives the binary cross-entropy objective~\cite{bishop2006pattern}. In case of the higher-order arrays, the low-rank representations are commonly constructed using the canonical polyadic (CP) decomposition, the Tucker decomposition, and the tensor-train (TT) decomposition~\cite{eckart1936,kolda2009,delathauwer2000,oseledets2011tt}. These formats generalize the idea of low-rank matrices through sums of the rank-one outer products, multilinear core-factor products, or tensor-network contractions. Their storage and approximation quality depend on the selected ranks. Reconstruction is commonly assessed with norm-based approximation errors and, for discrete or binary data, with entrywise and set-based measures such as Hamming error and intersection over union~\cite{kolda2009,sokolova2009metrics}. Although classical methods can represent binary data, discrete or nonlinear structures may require relatively high ranks to achieve a desired accuracy, motivating nonlinear low-rank representations.

Nonlinear matrix decomposition generalizes the classical model by applying an elementwise map to a low-rank product,
\[
\mathbf{A}\approx f(\mathbf{U}\mathbf{V}^{\mathsf T}).
\]
The nonlinear map changes the relation between the factors and the reconstructed entries and can therefore provide a more compact representation for data that are poorly matched by a direct linear approximation. Recently,  a latent-variable ReLU model for sparse nonnegative matrices has been introduced~\cite{saul2022}. This model establishes tractable inference, and has been derived via an expectation--maximization (EM) algorithm. Experiments on binary, grayscale, and color image data demonstrated substantially lower-rank representations than purely linear methods~\cite{saul2022}. Subsequent work developed accelerated schemes based on adaptive Nesterov extrapolation and three-block parameterizations~\cite{seraghiti2023}, momentum methods with regularization~\cite{wang2024}, and alternating partial Bregman methods for nonsmooth multiblock objectives~\cite{wang2025}. 

Later researchers separated the standard and latent ReLU formulations. They showed that block-coordinate and extrapolated algorithms converge under suitable assumptions~\cite{gillis2025}. Other nonlinear transformations include squaring and sum-of-squares constructions, while ADMM formulations separate the nonlinear map from the low-rank factors and allow different losses and clipping functions~\cite{awari2025admm}. For bounded data, sigmoid matrix decomposition has been introduced together with a block-coordinate optimization method~\cite{nguyen2025}. A subsequent formulation replaces the quadratic loss by binary cross-entropy. In such conditions, the resulting objective function is convex in one factor when the remaining factors are fixed~\cite{thayse2026}. Hence, the choice of the nonlinear map and discrepancy becomes a part of the decomposition model: it affects the geometry of the optimization problem, the interpretation of the reconstructed entries, and the convergence properties for the corresponding numerical algorithms.

Nonlinear tensor decomposition extends this construction to the higher-order arrays. A nonlinear CP model can be written as
\[
\widehat{\mathcal X}
=f\left(\sum_{q=1}^{r_{\mathrm{CP}}}\mathbf a_q^{(1)}\circ\mathbf a_q^{(2)}\circ\cdots\circ\mathbf a_q^{(d)}\right),
\]
where $r_{\mathrm{CP}}$ is the CP rank, $\mathbf a_q^{(k)}\in\mathbb{R}^{n_k}$ is the factor vector for component $q$ in mode $k$, $f$ acts elementwise, and $\circ$ denotes the outer product. The low-rank tensor is then an intermediate representation whose nonlinear image is compared with the observations. A ReLU–CP model for sparse nonnegative tensors was proposed and optimized by alternating minimization ~\cite{wang2026relu_ntd}. QR stabilization and dimension-tree contractions were used to control the numerical cost of tensor updates. 

A separate formulation considered clipped or saturated observations and introduced an auxiliary tensor together with an alternating direction method of multipliers (ADMM) scheme. In that setting, the nonlinear map models the observation mechanism, while the CP factors represent the underlying tensor rather than reproducing the clipped values directly~\cite{awari2026ncp}. These studies demonstrate how nonlinear tensor models can separate a low-rank multilinear representation from an entrywise observation transformation, but they focus on CP parameterizations with ReLU or clipping transformations. Tucker and TT formats provide structured parameterizations with explicit multilinear or train ranks~\cite{kolda2009,oseledets2011tt}. Motivated by these considerations, we extend nonlinear matrix decomposition to sigmoid Tucker and TT formats.

In this work, we develop sigmoid tensor decompositions for compression of binary data. Given a tensor $\mathcal B\in\{0,1\}^{n_1\times\cdots\times n_d}$, we represent a real-valued tensor $\mathcal Z$ in either Tucker or TT format and estimate its parameters by minimizing the binary cross-entropy between $\mathcal B$ and $\sigma(\mathcal Z)$. 
The resulting problem is jointly nonconvex in all factor matrices and cores, but it is convex in one active factor or core when the remaining blocks are fixed. We exploit this structure with block-coordinate gradient updates and form the binary reconstruction only after optimization by thresholding the sigmoid output. We compare the proposed models with linear Tucker and TT baselines and study the effects of tensor format, rank, initialization, and data geometry using different metrics. Our objective is to identify when sigmoid nonlinearization yields a compact and accurate representation of structured binary tensors and to provide a basis for further analysis of convergence, quantization, and nonlinear tensor parameterizations.

The remainder of the paper is organized as follows. Section~\ref{sec:problem}
introduces the notation and tensor formats and formulates the binary tensor
compression problem. Section~\ref{sec:models} introduces the sigmoid Tucker and
TT models, initialization, and block-coordinate optimization. Section~\ref{sec:protocol}
specifies the baselines, datasets, parameter budgets, and evaluation metrics.
Section~\ref{sec:results} presents the rank, initialization, and visual
reconstruction results. Section~\ref{sec:complexity} analyzes computational
cost. Section~\ref{sec:conclusion} concludes the paper. Additional visual
reconstructions are provided in Appendix~\ref{sec:appendix-visual}.

\section{Problem formulation}\label{sec:problem}

\subsection{Notation and tensor formats}\label{sec:notation-formats}

We use bold uppercase letters for matrices, calligraphic uppercase letters for
tensors, and lowercase letters for scalars. For a tensor
$\mathcal{X}\in\mathbb{R}^{n_1\times\cdots\times n_d}$, its entries are denoted
by $\mathcal{X}_{i_1,\ldots,i_d}$, where $1\leq i_m\leq n_m$.  For short notation,
we also write $\mathcal{X}_i=\mathcal{X}_{i_1,\ldots,i_d}$. Application of a scalar
function $f$ to tensor as $f(\mathcal{X})$ denotes its entrywise calculation. The mode-$m$ product of a tensor and a matrix is written as $\times_m$.

For a $d$-way tensor, $\operatorname{unfold}_m(\cdot)$ denotes the mode-$m$
matricization: the mode-$m$ index is placed along the rows and all remaining
indices are combined into columns. We use subscripts $(m)$ to denote these
unfoldings, for example $\mathbf{X}_{(m)}=\operatorname{unfold}_m(\mathcal X)$.
Throughout, $r=(r_1,\ldots,r_d)$ denotes the tuple of rank parameters. When
we refer to Tucker or TT rank $r$ without listing its components, all mode
ranks (or all internal TT ranks) are understood to have this common value.

The Tucker format represents a $d$-way tensor as
\begin{equation}
\mathcal{X}=\mathcal{G}\times_1\mathbf{U}^{(1)}
\times_2\mathbf{U}^{(2)}\cdots\times_d\mathbf{U}^{(d)},
\label{eq:tucker-format}
\end{equation}
where $\mathcal{G}\in\mathbb{R}^{r_1\times\cdots\times r_d}$ is the core and
$\mathbf{U}^{(m)}\in\mathbb{R}^{n_m\times r_m}$ are factor matrices satisfying
the orthogonality conditions
$\mathbf{U}^{(m)\mathsf T}\mathbf{U}^{(m)}=\mathbf{I}_{r_m}$, where
$\mathbf{I}_{r_m}$ is the $r_m\times r_m$ identity matrix. We use this
orthogonal Tucker parameterization in the algorithms below.
The Tucker ranks $(r_1,\ldots,r_d)$ determine the dimensions of the multilinear
subspaces. 

The tensor-train (TT) representation is given by
\begin{equation}
\mathcal{X}_{i_1,\ldots,i_d}=
\mathbf{G}^{(1)}[i_1]\mathbf{G}^{(2)}[i_2]\cdots\mathbf{G}^{(d)}[i_d],
\label{eq:tt-format}
\end{equation}
where the $k$th TT core is
$\mathcal{G}^{(k)}\in\mathbb{R}^{r_{k-1}\times n_k\times r_k}$ and
$\mathbf{G}^{(k)}[i_k]\in\mathbb{R}^{r_{k-1}\times r_k}$ is its slice at
index $i_k$. The boundary ranks satisfy $r_0=r_d=1$.

\subsection{Nonlinear tensor decomposition problem}\label{sec:general-problem}

Let $\mathcal{X}\in\mathbb{R}^{n_1\times\cdots\times n_d} $
be a $d$-way tensor and let $\mathfrak{M}_r$ denote a chosen low-rank tensor model with rank parameters collected in $r$. A nonlinear tensor decomposition seeks a real-valued representation tensor $\mathcal{Z}\in\mathfrak{M}_r$ such that
\begin{equation}
\min_{\mathcal{Z}\in\mathfrak{M}_r}
 D\bigl(\mathcal{X},f(\mathcal{Z})\bigr),
\label{eq:general-ntd}
\end{equation}
where $f$ is an entrywise nonlinear map and $D$ is a discrepancy function. The identity map gives the usual linear tensor approximation, while ReLU and sigmoid maps lead to nonlinear tensor models. In this work we focus on binary tensors $ \mathcal{B}\in\{0,1\}^{n_1\times\cdots\times n_d},     ~N_{\mathrm{tot}}=\prod_{\ell=1}^{d}n_\ell, $
and use a real-valued representation tensor $\mathcal{Z}$. Its entries are mapped to bounded values in $(0,1)$ through the sigmoid function
\[
\mathcal{P}=\sigma(\mathcal{Z}),
\qquad
\sigma(z)=\frac{1}{1+\exp(-z)}.
\]
Here $\mathcal{Z}$ is simply a real-valued representation of the binary tensor: it is the tensor before the sigmoid and thresholding operations, rather than an additional measured quantity. Large positive entries produce sigmoid values closer to one, while negative entries produce values closer to zero.
The sigmoid nonlinear tensor-decomposition problem is
\begin{equation}
\min_{\mathcal{Z}\in\mathfrak{M}_r}
\mathcal{L}_{\mathrm{CE}}\bigl(\mathcal{B},\sigma(\mathcal{Z})\bigr),
\label{eq:sigmoid-ntd}
\end{equation}
where $\mathcal{Z}$ has the selected Tucker or TT representation. 

The objective loss function is jointly nonconvex in all factor
blocks; with all blocks except one factor matrix, Tucker core, or TT core fixed,
it is convex in the active block. The average binary cross-entropy loss is
\begin{equation}
\mathcal{L}_{\mathrm{CE}}(\mathcal{B},\sigma(\mathcal{Z}))
=-\frac{1}{N_{\mathrm{tot}}}\sum_{i_1=1}^{n_1}\cdots\sum_{i_d=1}^{n_d}
\left[
\mathcal{B}_{i_1,\ldots,i_d}\log\sigma(\mathcal{Z}_{i_1,\ldots,i_d})
+\bigl(1-\mathcal{B}_{i_1\ldots i_d}\bigr)
\log\bigl(1-\sigma(\mathcal{Z}_{i_1,\ldots,i_d})\bigr)
\right].
\label{eq:binary-ce}
\end{equation}
Equivalently, one may also write it as 
\vspace{-0.5cm}
\[
\mathcal{L}_{\mathrm{CE}}
=\frac{1}{N_{\mathrm{tot}}}\sum_{i_1=1}^{n_1}\cdots\sum_{i_d=1}^{n_d}
\left[
\log\bigl(1+\exp(\mathcal{Z}_{i_1,\ldots,i_d})\bigr)
-\mathcal{B}_{i_1,\ldots,i_d}\mathcal{Z}_{i_1,\ldots,i_d}
\right].
\]
The decoded binary tensor is determined by the sigmoid threshold $1/2$, or equivalently by the sign of the representation values:
\begin{equation}
\widehat{\mathcal{B}}_{i_1,\ldots,i_d}
=\mathbf{1}\left\{\mathcal{P}_{i_1,\ldots,i_d}\geq\frac12\right\}
=\mathbf{1}\left\{\mathcal{Z}_{i_1,\ldots,i_d}\geq0\right\},
\label{eq:binary-decoder}
\end{equation}
where $\mathbf{1}\{\cdot\}$ denotes the indicator function.




\section{Nonlinear tensor decomposition}\label{sec:models}

In this section, we describe the sigmoid
Tucker and sigmoid tensor-train models for binary-data compression.

\subsection{Sigmoid Tucker model}
Using the Tucker factors and core defined in Section~\ref{sec:notation-formats},
the real-valued representation is
\begin{equation}
\mathcal{Z}=\mathcal{G}\times_1\mathbf{U}^{(1)}\times_2\cdots\times_d\mathbf{U}^{(d)}.
\label{eq:sigmoid-tucker}
\end{equation}
The reconstructed tensor is obtained from $\sigma(\mathcal{Z})$ and thresholded according to~\eqref{eq:binary-decoder}. The sigmoid Tucker model minimizes the cross-entropy objective in~\eqref{eq:binary-ce} with respect to the core and factor matrices; its implementation is summarized in Algorithm~\ref{alg:ce-tucker}.
\begin{algorithm}[H]
\caption{Sigmoid Tucker compression (CE--Tucker)}
\label{alg:ce-tucker}
\begin{algorithmic}[1]
\Require Binary tensor $\mathcal{B}$; Tucker rank $\boldsymbol r$; maximum number of iterations $N_{\mathrm{iter}}$.
\Ensure Tucker factors $\{\mathbf U^{(m)}\}_{m=1}^{d}$ and core $\mathcal G$.
\State Initialize $\mathcal G,\{\mathbf U^{(m)}\}_{m=1}^{d}$ using the selected strategy. \Comment{Section~\ref{sec:results-initialization}.}
\State Orthogonalize the initial factor matrices and absorb the transformations into $\mathcal G$.
\State Compute the current representation $\mathcal Z\gets\mathcal G\times_1\mathbf U^{(1)}\cdots\times_d\mathbf U^{(d)}$.
\For{$k=1,\ldots,N_{\mathrm{iter}}$}
    \For{$m=1,\ldots,d$}
        \State $\mathcal Z\gets\mathcal G\times_1\mathbf U^{(1)}\cdots\times_d\mathbf U^{(d)}$
        \State  $\mathcal E\gets\sigma(\mathcal Z)-\mathcal B$ and
        \State $\mathbf K_m\gets\mathbf U^{(d)}\otimes\cdots\otimes\mathbf U^{(m+1)}\otimes\mathbf U^{(m-1)}\otimes\cdots\otimes\mathbf U^{(1)}$.
        \State $\mathbf E_{(m)}\gets\operatorname{unfold}_m(\mathcal E)$
        \State $\mathbf G_{(m)}\gets\operatorname{unfold}_m(\mathcal G)$.
        \State $\mathbf U^{(m)}\gets\mathbf U^{(m)}-\dfrac{\nu_U}{N_{\mathrm{tot}}}\mathbf E_{(m)}\mathbf K_m\mathbf G_{(m)}^{\mathsf T}$.
    \EndFor
    \State Orthogonalize the updated factor matrices and absorb the transformations into $\mathcal G$.
    \State $\mathcal Z\gets\mathcal G\times_1\mathbf U^{(1)}\cdots\times_d\mathbf U^{(d)}$
    \State $\mathcal E\gets\sigma(\mathcal Z)-\mathcal B$.
    \State $\mathcal G\gets\mathcal G-\dfrac{\nu_G}{N_{\mathrm{tot}}}\mathcal E\times_1\mathbf U^{(1)\mathsf T}\cdots\times_d\mathbf U^{(d)\mathsf T}$.
\EndFor
\State Form the final representation $\mathcal Z\gets\mathcal G\times_1\mathbf U^{(1)}\cdots\times_d\mathbf U^{(d)}$.
\State Compute $\mathcal P\gets\sigma(\mathcal Z)$ and $\widehat{\mathcal B}_{i_1,\ldots,i_d}\gets\mathbf 1\{\mathcal P_{i_1,\ldots,i_d}\geq1/2\}$.
\State \Return $\mathcal G,\{\mathbf U^{(m)}\}_{m=1}^{d}$ and $\widehat{\mathcal B}$.
\end{algorithmic}
\end{algorithm}

The step sizes are chosen from curvature bounds for the active Tucker block
rather than tuned separately for each object or rank, no momentum or adaptive
step-size search is used. 

This block-coordinate strategy follows the general
alternating optimization framework used in nonlinear matrix decomposition
~\cite{gillis2025}. The factor update in Algorithm~\ref{alg:ce-tucker} is
written with a factor $1/N_{\mathrm{tot}}$ because the objective in~\eqref{eq:binary-ce} is
the average cross-entropy over the $N_{\mathrm{tot}}=\prod_{m=1}^d n_m$ tensor entries. Thus,
if the gradient contraction is formed as a sum over entries, division by $N_{\mathrm{tot}}$
converts it to the gradient of the averaged objective. The remaining step-size
constants are selected from an upper bound on the corresponding block Hessian
(the sigmoid derivative contributes at most $1/4$), together with the norm of
the fixed multilinear design matrix.

After a Tucker factor update, QR orthogonalization removes the scaling
ambiguity between the factor and core and improves conditioning of the next
block subproblem, as in orthogonal Tucker/HOSVD methods~\cite{delathauwer2000}.
The triangular factor is absorbed into the core, so the represented tensor is
unchanged. The core itself is not orthogonalized: it stores the multilinear
coefficients and is updated by its own cross-entropy gradient step. We denote
the corresponding step sizes by $\nu_U$ for the Tucker factors and $\nu_G$ for
the Tucker core.

\subsection{Sigmoid tensor-train model}
For the TT cores and ranks defined in Section~\ref{sec:notation-formats}, the
representation entries are
\begin{equation}
\mathcal{Z}_{i_1,\ldots,i_d}
=\mathbf{G}^{(1)}[i_1]\cdots\mathbf{G}^{(d)}[i_d].
\label{eq:sigmoid-tt}
\end{equation}
The sigmoid TT model uses the same sigmoid decoder and cross-entropy objective as the Tucker model.
The corresponding optimization procedure is given in Algorithm~\ref{alg:ce-tt}.

For clarity, before updating the active core $\mathbf G^{(m)}$, we use the
standard TT gauge transformation to put the cores into mixed-canonical form
with the orthogonality center at position $m$. The cores to the left satisfy
the left-orthogonality conditions
\begin{equation}
\sum_{i_j=1}^{n_j}\mathbf G^{(j)}[i_j]^{\mathsf T}
\mathbf G^{(j)}[i_j]=\mathbf I_{r_j},
\qquad j<m,
\end{equation}
whereas the cores to the right satisfy
\begin{equation}
\sum_{i_j=1}^{n_j}\mathbf G^{(j)}[i_j]
\mathbf G^{(j)}[i_j]^{\mathsf T}=\mathbf I_{r_{j-1}},
\qquad j>m.
\end{equation}
These conditions are imposed by QR factorizations again. For a left core, we reshape
the slices into a matrix of size $(r_{j-1}n_j)\times r_j$ and compute
$\mathbf A_j=\mathbf Q_j\mathbf R_j$; the core is replaced by the reshaped
$\mathbf Q_j$, and $\mathbf R_j$ is
absorbed into the next core. For a right core, we factor the transpose of the
reshaped matrix, $\mathbf A_j^{\mathsf T}=\mathbf Q_j\mathbf R_j$, replace the
core by the reshaped $\mathbf Q_j^{\mathsf T}$, and absorb
$\mathbf R_j^{\mathsf T}$ into the preceding core.
This trick changes only the internal TT gauge and leaves
\(\mathcal Z_{i_1,\ldots,i_d}\) unchanged. Such type of orthogonalization is
standard in TT-alternating optimization and local core-update schemes
~\cite{oseledets2011tt,holtz2012als,rohwedder2013als}. It improves conditioning of the active-core gradient by removing arbitrary
scalings from the neighboring cores.

\begin{algorithm}[H]
\caption{Sigmoid tensor-train compression (CE--TT)}
\label{alg:ce-tt}
\begin{algorithmic}[1]
\Require Binary tensor $\mathcal{B}$; TT rank $r$; maximum number of iterations $N_{\mathrm{iter}}$.
\Ensure TT cores $\{\mathbf G^{(m)}\}_{m=1}^{d}$.
\State Initialize $\{\mathbf G^{(m)}\}_{m=1}^{d}$ using the selected strategy. \Comment{Section~\ref{sec:results-initialization}.}

\For{$k=1,\ldots,N_{\mathrm{iter}}$}
    \For{$m=1,\ldots,d$}
        \State Put TT cores into mixed-canonical form with orthogonality center at $m$.
        \State $\mathcal Z_{i_1,\ldots,i_d}\gets\mathbf G^{(1)}[i_1]\cdots\mathbf G^{(d)}[i_d]$
        \State  $\mathcal E\gets\sigma(\mathcal Z)-\mathcal B$.
        \State $\boldsymbol\ell_{i_{<m}}\gets\mathbf G^{(1)}[i_1]\cdots\mathbf G^{(m-1)}[i_{m-1}] \forall\,i_m:$
        \State $\boldsymbol\rho_{i_{>m}}\gets\mathbf G^{(m+1)}[i_{m+1}]\cdots\mathbf G^{(d)}[i_d]\qquad \forall i_m:$
        \State $\mathbf G^{(m)}[i_m]\gets\mathbf G^{(m)}[i_m]-\dfrac{\nu_m}{N_{\mathrm{tot}}}\sum_{i_{\ne m}}\mathcal E_{i_1\ldots i_d}\,\boldsymbol\ell_{i_{<m}}^{\mathsf T}\boldsymbol\rho_{i_{>m}}^{\mathsf T}$.
    \EndFor
\EndFor
\State Form the final representation $\mathcal Z_{i_1,\ldots,i_d}\gets\mathbf G^{(1)}[i_1]\cdots\mathbf G^{(d)}[i_d]$.
\State Compute $\mathcal P\gets\sigma(\mathcal Z)$ and $\widehat{\mathcal B}_{i_1,\ldots,i_d}\gets\mathbf 1\{\mathcal P_{i_1,\ldots,i_d}\geq1/2\}$.
\State \Return $\{\mathbf G^{(m)}\}_{m=1}^{d}$ and $\widehat{\mathcal B}$.
\end{algorithmic}
\end{algorithm}

\subsection{Initialization}
The initialization determines the starting point of the gradient-based block updates and can therefore affect both the final loss and the reconstruction quality. Similar finite-logit and decomposition-based starts were used for sigmoid matrix decomposition with cross-entropy~\cite{thayse2026}. 

Our main initialization replaces each binary entry by a finite logit corresponding to a sigmoid output of $0.98$ for one and $0.02$ for zero. Then, we apply the selected linear tensor decomposition to this real-valued tensor and use its factors as the initial parameters. This procedure provides an informative starting sign pattern while retaining a low-rank parameterization. The effects of these initialization choices on convergence and reconstruction quality are studied in Section~\ref{sec:results-initialization}.

\section{Experimental protocol}\label{sec:protocol}

This section specifies the datasets, comparison methods, optimization parameters, and evaluation criteria used to assess the sigmoid tensor decompositions. The main rank sweep fixes the same numerical rank across Tucker and TT and reports their different parameter counts explicitly.

\subsection{Datasets}

\begin{enumerate}
  \item \textbf{ModelNet40 objects.} We use 17 voxelized CAD objects from five categories: eight chairs, six cars, one airplane, one sofa, and one table. The original CAD collection is available from the \href{https://modelnet.cs.princeton.edu/}{Princeton ModelNet40 dataset}. Each object is represented by a binary occupancy tensor of size $128\times128\times128$. The occupancy tensors are obtained from the corresponding signed-distance representation using
  \[
  \mathcal{B}_{i,j,k}=\mathbf{1}\{\mathcal{X}_{i,j,k}<0\}.
  \]
  \item \textbf{Stanford scans.} We use the Dragon and Armadillo scans from the \href{http://graphics.stanford.edu/data/3Dscanrep/}{Stanford 3D Scanning Repository} as a separate geometric group. The source scans are converted to binary occupancy tensors and evaluated at size $128^3$ in the aggregate rank study. They are not pooled with ModelNet40 when computing dataset means.
  \item \textbf{Medical Segmentation Decathlon masks.} We use ten Heart masks (Task02), ten Hippocampus masks (Task04), and ten Spleen masks (Task09) from the official \href{https://medicaldecathlon.com/}{Medical Segmentation Decathlon}. Native mask dimensions vary between cases. Each foreground label is converted to a binary mask and resampled to $128^3$ with nearest-neighbor interpolation. The Hippocampus anterior and posterior labels are merged. These experiments compress segmentation masks rather than the original medical-image intensity volumes.
\end{enumerate}

The aggregate rank experiments therefore contain three data groups: ModelNet40, Stanford, and medical masks. The high-resolution visualization study uses selected Chair and Dragon tensors at $512^3$; those visualizations are reported separately from the $128^3$ aggregate rank sweep.

\subsection{Methods}

The comparison uses two proposed sigmoid models and three direct tensor baselines
described in Table~\ref{tab:methods-protocol}. The baselines isolate the effect of
sigmoid nonlinearization from the effect of the tensor format.

\begin{table}[h]
\centering
\caption{Methods used in the binary-tensor comparison.}
\label{tab:methods-protocol}
\small
\begin{tabular}{p{0.22\textwidth}p{0.64\textwidth}}
\toprule
Method & Representation and objective\\
\midrule
TT-TSDF & TT-SVD applied to a continuous signed-distance field; decode by its zero level set ($\widehat{\mathcal B}=\mathbf{1}\{\widehat{\mathcal X}<0\}$), following the continuous-field compression strategy of Boyko et al.~\cite{boyko2020}.\\
TT-binary & TT-SVD applied directly to the binary tensor.\\
Tucker-binary &  ST-HOSVD applied directly to the binary tensor.\\
CE-Tucker & Sigmoid Tucker optimized with binary cross-entropy.\\
CE-TT & Sigmoid tensor-train optimized with binary cross-entropy.\\
\bottomrule
\end{tabular}
\end{table}

For TT-binary and Tucker-binary, the decompressed entries are real-valued
approximations to the stored binary entries, so we obtain the binary output by
thresholding at $0.5$. In contrast, TT-TSDF reconstructs a signed-distance field;
its occupied region is defined by the zero level set, with negative signed
distance values classified as foreground.

Tucker and TT are standard multilinear formats~\cite{kolda2009,delathauwer2000,oseledets2011tt}. ST-HOSVD and TT-SVD provide deterministic starts for the direct baselines and the default structured starts for the sigmoid models. Alternative initializations are examined in the initialization study.

For the rank sweep, $r$ denotes the common numerical rank used for both formats. Their storage costs are different:
\[
 p_{\mathrm{TT}}(r)=n(r^2+2r), \qquad
 p_{\mathrm{Tucker}}(r)=r^3+3nr, \qquad r\in\{4,8,12,16\}.
\]
For a third-order tensor with equal mode sizes $n_1=n_2=n_3=n$, the formulas
above give the corresponding three-dimensional storage budgets. For Tucker, $r$ means $r_1=r_2=r_3=r$; for TT, it means the two internal ranks satisfy
$r_1=r_2=r$ with boundary ranks $r_0=r_3=1$.
Thus, equal rank does not imply equal storage. The actual parameter count is reported for every result, and the parameter-error figure uses this count on its horizontal axis. All factors are stored in FP32 for the reported storage comparison.

Unless otherwise stated, nonlinear fits use at most 1000 block-coordinate sweeps and record the optimization history every 20 sweeps. Reaching the iteration limit is treated as a finite-budget result rather than a global convergence certificate.

\subsection{Metrics}

Let $N_{\mathrm{tot}}$ be the number of voxels, let $B_i\in\{0,1\}$ be the reference value, let $Z_i$ be the reconstructed real-valued representation, and let $\widehat{B}_i$ be the decoded binary value. We report the following metrics.
Their definitions and interpretations are summarized in Table~\ref{tab:metrics}.

\begin{table}[h]
\centering
\caption{Metrics used in the binary-tensor experiments.}
\label{tab:metrics}
\small
\begin{tabular}{p{0.18\textwidth}p{0.38\textwidth}p{0.35\textwidth}}
\toprule
Metric & Mathematical definition & Interpretation \\
\midrule
Binary cross-entropy (CE) & $\displaystyle \frac{1}{N_{\mathrm{tot}}}\sum_{i=1}^{N_{\mathrm{tot}}}\log\!\left(1+\exp\!\left(-(2B_i-1)Z_i\right)\right)$ & Objective optimized by the sigmoid models; lower values indicate better agreement with the binary targets. \\
Normalized Frobenius error & $\displaystyle e_F=\frac{\|\widehat{\mathcal{B}}-\mathcal{B}\|_F}{\|\mathcal{B}\|_F}$ & Relative error of the decoded tensor. \\
Hamming error & $\displaystyle e_H=\frac{1}{N_{\mathrm{tot}}}\sum_{i=1}^{N_{\mathrm{tot}}}\mathbf{1}\{\widehat{B}_i\ne B_i\}$ & Fraction of voxels whose decoded binary value is incorrect. \\
Intersection over union (IoU) & $\displaystyle \operatorname{IoU}=\frac{|\widehat{\mathcal{B}}\cap\mathcal{B}|}{|\widehat{\mathcal{B}}\cup\mathcal{B}|}$ & Overlap between reconstructed and reference occupied regions; higher values indicate better shape agreement. \\
\bottomrule
\end{tabular}
\end{table}

IoU, Hamming error, and normalized Frobenius error are defined for all five comparison methods (see Table \ref{tab:metrics}). Cross-entropy is reported for the sigmoid models when a continuous score tensor is available. We report object-level or case-level means and standard deviations within each dataset group, together with the number of stored parameters.

\section{Results}\label{sec:results}

The results are organized around three questions: 
\begin{enumerate}
    \item how rank affect reconstruction, 
    \item how initialization affects the nonconvex optimization, \item how the numerical metrics appear in three-dimensional reconstructions.
\end{enumerate}

The initialization study uses representative Chair, Dragon, and Heart cases. For the main visualization, we show the ModelNet40 Chair at rank $r=8$; the remaining visual examples are moved to Appendix~\ref{sec:appendix-visual}.

\subsection{Rank and storage-budget sensitivity}\label{sec:results-rank}

The aggregate rank study contains 17 ModelNet40 objects, two Stanford scans, and 30 medical masks from the Heart, Hippocampus, and Spleen tasks. Each tensor in this study has size $128^3$. We evaluate $r\in\{4,8,12,16\}$ for the five methods in Section~\ref{sec:protocol}. The numerical rank is shared between Tucker and TT, but their parameter counts differ and are reported through the storage-error plots.

Table~\ref{tab:rank-sweep-iou} reports mean intersection over union (IoU) and standard deviation. The sigmoid models improve strongly as the rank increases. CE-TT is the strongest method in the aggregate IoU comparison, while CE-Tucker also improves consistently and becomes highly accurate at moderate rank. The direct binary baselines improve more gradually. The same trend is visible in the normalized Frobenius and Hamming errors.

\begin{figure}[!ht]
\centering
\includegraphics[width=\textwidth]{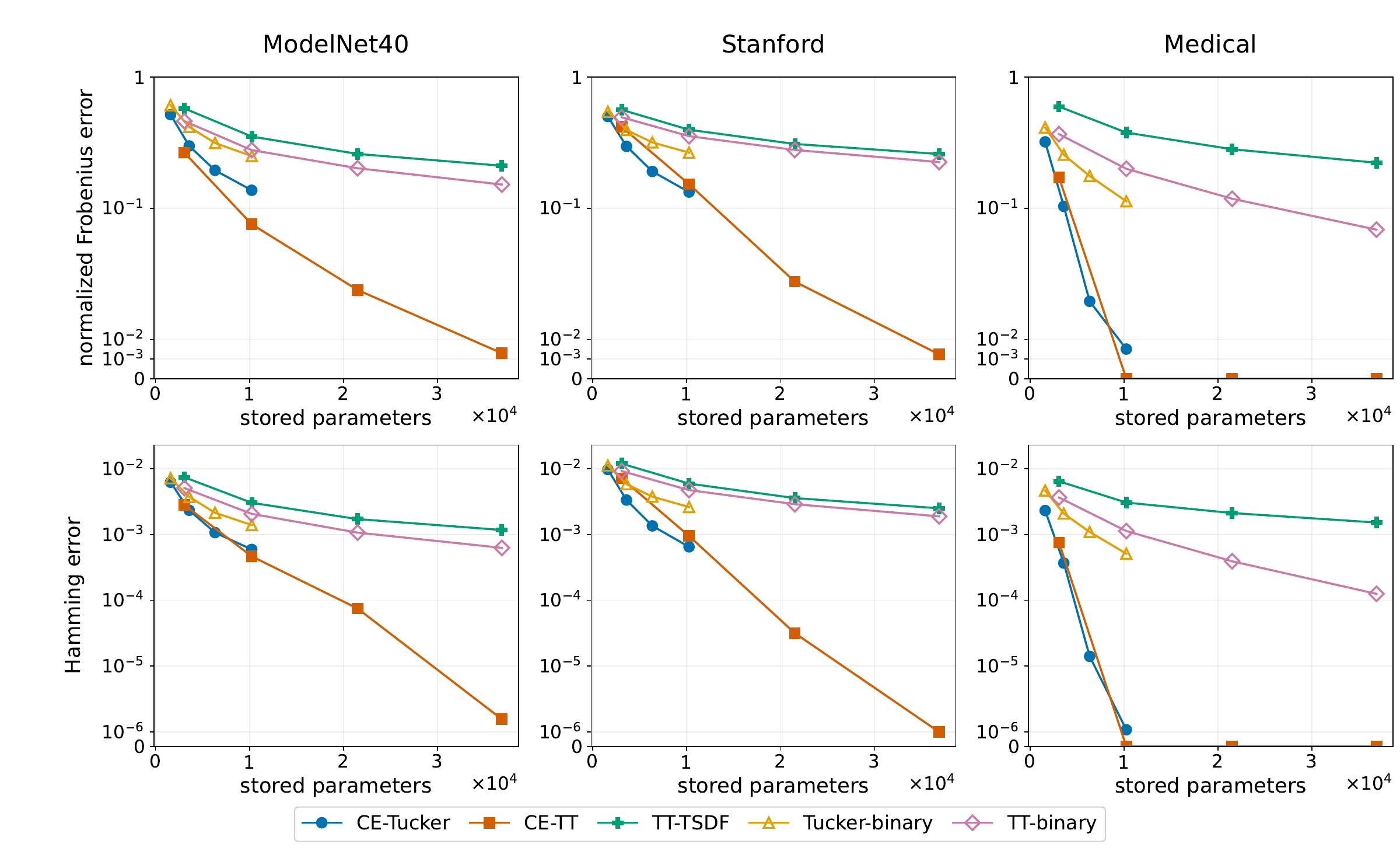}
\caption{Normalized Frobenius and Hamming error versus stored parameters. Rows show the two errors and columns show ModelNet40, Stanford, and pooled medical masks. Each curve contains $r\in\{4,8,12,16\}$.}
\label{fig:rank-sweep-all}
\end{figure}

The stronger aggregate performance of CE-TT is partly explained by its larger
number of free parameters at the same numerical rank. This difference is visible
in the storage-error plots in Figure~\ref{fig:rank-sweep-all}. We see that equal Tucker and TT
ranks do not correspond to equal storage budgets. In these rank-sweep experiments,
CE-Tucker is initialized with ST-HOSVD applied to the finite-logit target, whereas
CE-TT is initialized with the corresponding TT-SVD decomposition. The separate
initialization study evaluates how this choice affects convergence. The decreasing
standard deviations as the rank increases indicate that the objects become more
uniformly well represented at higher ranks. This behavior is consistent with the
datasets possessing an effective low-rank structure: increasing the rank moves the
models toward the intrinsic ranks of the data and reduces variation in reconstruction
quality across objects.

\begin{table}[h]
\centering
\caption{Mean IoU $\pm$ standard deviation for the binary-tensor rank sweep. The columns use the common numerical rank $r$; Tucker and TT have different storage costs at each rank.}
\label{tab:rank-sweep-iou}
\scriptsize
\newcommand{\datasetseparator}{\multicolumn{6}{c}{\makebox[0.98\linewidth]{\leaders\hbox to 6pt{\hss\textendash\hss}\hfill}}\\[-2pt]}
\begin{tabular}{llrrrr}\toprule
Dataset & Method & $r=4$ & $r=8$ & $r=12$ & $r=16$\\\midrule
ModelNet40 & CE-Tucker & 0.720 $\pm$ 0.182 & 0.892 $\pm$ 0.094 & 0.951 $\pm$ 0.047 & 0.974 $\pm$ 0.027 \\
ModelNet40 & CE-TT & \textbf{0.861 $\pm$ 0.128} & \textbf{0.976 $\pm$ 0.037} & \textbf{0.996 $\pm$ 0.010} & \textbf{1.000 $\pm$ 0.000} \\
ModelNet40 & Tucker-binary & 0.637 $\pm$ 0.167 & 0.814 $\pm$ 0.127 & 0.888 $\pm$ 0.088 & 0.926 $\pm$ 0.066 \\
ModelNet40 & TT-binary & 0.690 $\pm$ 0.148 & 0.856 $\pm$ 0.108 & 0.923 $\pm$ 0.071 & 0.953 $\pm$ 0.048 \\
ModelNet40 & TT-TSDF & 0.526 $\pm$ 0.184 & 0.788 $\pm$ 0.142 & 0.882 $\pm$ 0.091 & 0.920 $\pm$ 0.065 \\
\datasetseparator
Stanford & CE-Tucker & 0.766 $\pm$ 0.092 & 0.914 $\pm$ 0.024 & 0.964 $\pm$ 0.004 & 0.982 $\pm$ 0.000 \\
Stanford & CE-TT & \textbf{0.831 $\pm$ 0.108} & \textbf{0.976 $\pm$ 0.016} & \textbf{0.999 $\pm$ 0.001} & \textbf{1.000 $\pm$ 0.000} \\
Stanford & Tucker-binary & 0.724 $\pm$ 0.052 & 0.849 $\pm$ 0.024 & 0.902 $\pm$ 0.012 & 0.931 $\pm$ 0.004 \\
Stanford & TT-binary & 0.769 $\pm$ 0.051 & 0.878 $\pm$ 0.026 & 0.924 $\pm$ 0.012 & 0.950 $\pm$ 0.007 \\
Stanford & TT-TSDF & 0.700 $\pm$ 0.055 & 0.845 $\pm$ 0.035 & 0.905 $\pm$ 0.016 & 0.933 $\pm$ 0.009 \\
\datasetseparator
Medical & CE-Tucker & 0.897 $\pm$ 0.038 & 0.988 $\pm$ 0.009 & 0.999 $\pm$ 0.001 & \textbf{1.000 $\pm$ 0.000} \\
Medical & CE-TT & \textbf{0.967 $\pm$ 0.020} & \textbf{1.000 $\pm$ 0.000} & \textbf{1.000 $\pm$ 0.000} & \textbf{1.000 $\pm$ 0.000} \\
Medical & Tucker-binary & 0.834 $\pm$ 0.046 & 0.933 $\pm$ 0.030 & 0.967 $\pm$ 0.019 & 0.985 $\pm$ 0.010 \\
Medical & TT-binary & 0.868 $\pm$ 0.039 & 0.959 $\pm$ 0.017 & 0.985 $\pm$ 0.008 & 0.995 $\pm$ 0.004 \\
Medical & TT-TSDF & 0.669 $\pm$ 0.090 & 0.865 $\pm$ 0.024 & 0.922 $\pm$ 0.020 & 0.950 $\pm$ 0.020 \\
\bottomrule
\end{tabular}

\end{table}

Figure~\ref{fig:rank-sweep-all} shows the normalized Frobenius and Hamming errors against the actual number of stored parameters. Thus, the figure compares the methods at their true storage costs rather than treating equal numerical ranks as equal budgets. The three columns correspond to ModelNet40, Stanford, and the pooled medical masks.

The rank study separates two effects. Increasing $r$ enlarges the representable tensor class, while changing the format changes the storage budget and contraction structure. A low IoU at small rank therefore indicates a capacity limit or a poor fit, whereas a high Hamming error indicates that many individual voxels are misclassified. Since the data are sparse, IoU is especially informative for thin structures: a small number of missed occupied voxels can produce a large overlap loss even when the global Hamming error remains small.

\subsection{Initialization sensitivity and convergence}\label{sec:results-initialization}

This subsection examines how the
optimization behaves under different initialization strategies. We compare random,
decomposition, rank-expansion warm start, and residual-informed initialization. The decomposition
start uses ST-HOSVD for Tucker and TT-SVD for TT applied to the finite-logit target.
Rank-expansion warm start expands a fitted lower-rank model, while residual-informed
initialization compresses the unresolved finite-logit residual before adding the new
rank blocks. Their definitions are summarized in Table~\ref{tab:initialization-methods}.

\begin{figure}[!ht]
\centering
\includegraphics[width=\textwidth]{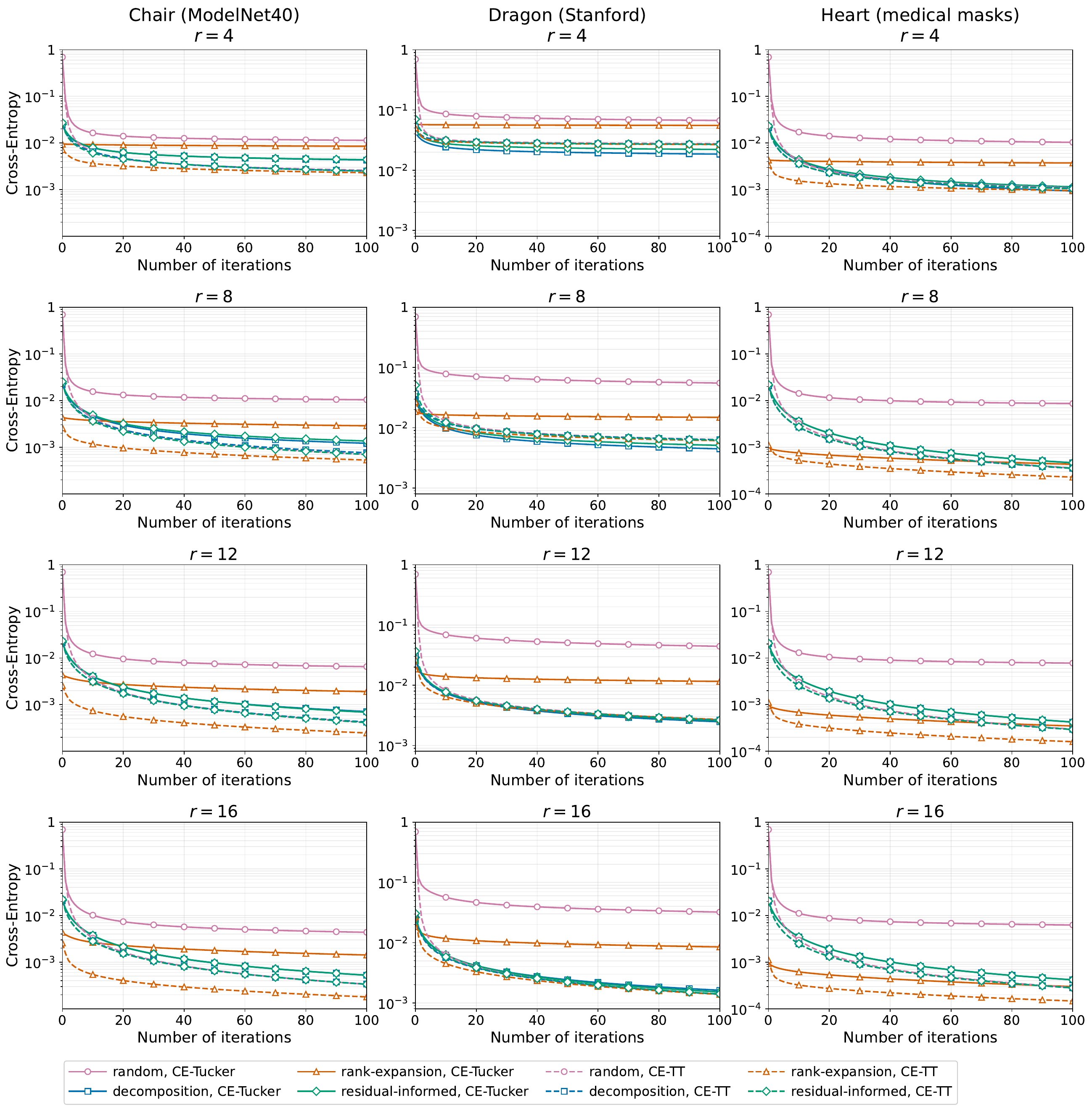}
\caption{Initialization sensitivity for CE-Tucker and CE-TT. Rows correspond to $r=4,8,12,16$ and columns correspond to Chair, Dragon, and Heart. The curves show cross-entropy over 100 block-coordinate iterations.}
\label{fig:initialization-grid-four-ranks}
\end{figure}

\begin{table}[h]
\centering
\caption{Initialization strategies used for sigmoid Tucker and TT optimization.}
\label{tab:initialization-methods}
\small
\begin{tabular}{p{0.20\linewidth}p{0.70\linewidth}}
\toprule
Strategy & Construction \\
\midrule
Random & Random factor blocks and a scaled random core at the target rank. \\
Decomposition & ST-HOSVD or TT-SVD of the finite-logit target, with the resulting factors and core used directly as the initial parameters. \\
Rank-expansion warm start & First construct a lower-rank decomposition with rank $r_0<r$, then expand its factors and core to the target rank by appending initialized components while preserving the lower-rank representation. Rank-adaptive expansion of this type is used in tensor-train optimization~\cite{holtz2012als}. \\
Residual-informed & Construct a lower-rank approximation, compute the unresolved finite-logit residual, decompose that residual, and insert the resulting directions into the newly added factor blocks or TT cores. Residual-based local enrichment is standard in adaptive TT methods~\cite{holtz2012als,DolgovSavostyanov2014}. \\
\bottomrule
\end{tabular}
\end{table}
\vspace{-1cm}
Figure~\ref{fig:initialization-grid-four-ranks} reports the cross-entropy histories for three representative objects---Chair, Dragon, and Heart---at four ranks, $r\in\{4,8,12,16\}$. The results show that initialization interacts with both the tensor format and the object geometry.

For the Chair, CE-Tucker with random initialization gives the weakest convergence, and this behavior changes little as the rank increases. One possible explanation is that random Tucker factors select directions that have little overlap with the informative low-rank subspaces, so the corresponding block gradients provide limited progress. This explanation is incomplete, however, because CE-TT converges reasonably from the same type of random start. The difference may therefore also involve the conditioning and gauge structure of the two parameterizations. A precise explanation requires a separate analysis of the block Hessians and is left for future work. For CE-TT, the rank-expansion warm start is the most effective initialization on the Chair: the lower-rank fit supplies a useful geometric approximation, and the added components can refine it progressively. Decomposition and residual-informed starts are more stable across ranks for both formats.

For the Dragon, most initialization curves are close across the four ranks, with the main deviations occurring for two CE-Tucker starts. This indicates that the Dragon geometry is already well aligned with several of the tested initial subspaces, while the Tucker parameterization remains sensitive to particular random or rank-expansion choices. The TT curves are comparatively consistent, and no single initialization dominates as strongly as it does for the Chair.
The Heart results resemble the Chair results in that CE-Tucker is more initialization-sensitive. In this case, however, CE-Tucker with the rank-expansion warm start converges more favorably. The improvement is consistent with the medical masks having a smaller effective rank: after approximately rank eight, the dominant structure is already captured, so expanding a lower-rank solution provides a suitable starting point and leaves only a small residual to optimize. Overall, Figure~\ref{fig:initialization-grid-four-ranks} shows that convergence behavior depends jointly on the tensor format, the object geometry, and the initialization strategy. A mathematical analysis of this initialization-dependent convergence is left for future work.

The decomposition-based start is the most stable and reproducible initialization
across the three objects and the tested ranks. It can also be computed more
efficiently by replacing the deterministic ST-HOSVD and TT-SVD steps with
randomized range-SVD constructions using oversampling and one power iteration,
following the randomized tensor
approximation framework of Che and Wei~\cite{chewei2019randomized}. Table~\ref{tab:randomized-initialization}
reports the results averaged over the 35 available ModelNet objects. Each
randomized entry averages five independent seeds, whereas the deterministic
decomposition is deterministic for a fixed input. All initialization timings were measured 
on a single NVIDIA A100 GPU using PyTorch with CUDA. The reported time includes 
only the construction of the ST-HOSVD, TT-SVD, randomized ST-HOSVD, or randomized TT-SVD 
initialization and excludes data loading and subsequent cross-entropy optimization. The nonlinear optimization
uses the same 100 block-coordinate sweeps in all cases. Randomization changes
the final CE and IoU only in the third or fourth decimal place, while reducing
the initialization time substantially. Thus, when the ranks are prescribed in
advance, randomized SVD provides a cheaper structure-aligned initialization
without a measurable loss of reconstruction quality.

\begin{table}[!h]
\centering
\caption{Deterministic and randomized decomposition initialization on 35 ModelNet objects. }
\label{tab:randomized-initialization}
\small
\begin{tabular}{llrccc}
\toprule
Format & Rank & Initialization & Final CE & Final IoU & Init. time (s) \\
\midrule
Tucker & 8  & ST-HOSVD & 0.00570 & 0.8871 & 0.0370 \\
        &    & Randomized ST-HOSVD & 0.00571 & 0.8874 & 0.0051 \\
\midrule
Tucker & 16 & ST-HOSVD & 0.00197 & 0.9708 & 0.0314 \\
        &    & Randomized ST-HOSVD & 0.00198 & 0.9711 & 0.0063 \\
\midrule
TT & 8  & TT-SVD & 0.00196 & 0.9767 & 0.0297 \\
   &    & Randomized TT-SVD & 0.00196 & 0.9770 & 0.0041 \\
\midrule
TT & 16 & TT-SVD & 0.00047 & 0.9991 & 0.0623 \\
   &    & Randomized TT-SVD & 0.00047 & 0.9991 & 0.0046 \\
\bottomrule
\end{tabular}
\end{table}

\subsection{Visual reconstruction at the selected rank}\label{sec:results-visual}

For the main visual comparison, we select the ModelNet40 Chair at $512^3$ and $r=8$. The figure includes the reference and direct baselines in the first row, followed by CE-Tucker and CE-TT reconstructions for the four initialization strategies. The IoU value is shown for each reconstruction.

\begin{figure}[!ht]
\centering
\includegraphics[width=\textwidth]{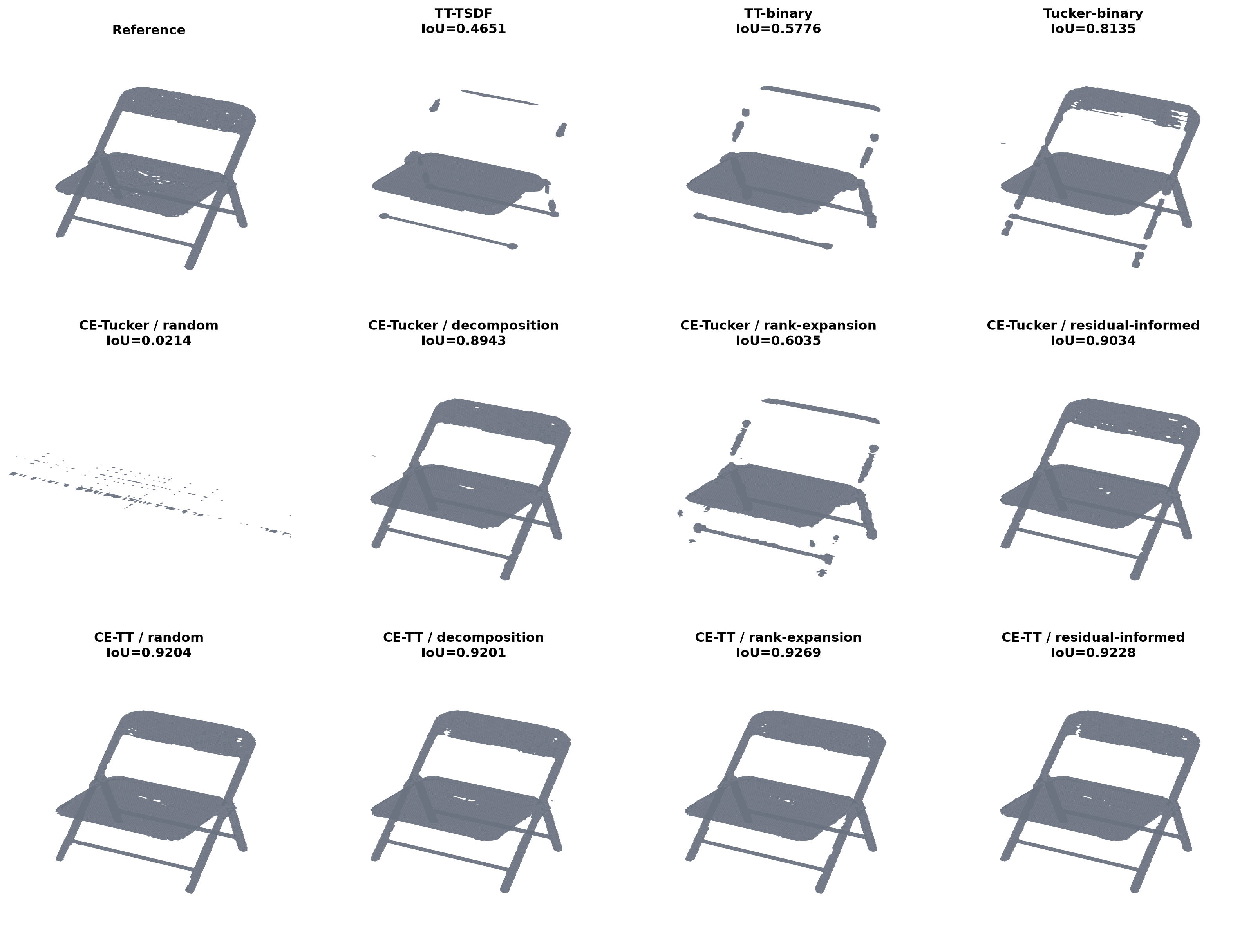}
\caption{Three-dimensional reconstruction of the $512^3$ ModelNet40 Chair at rank $r=8$.}
\label{fig:final-chair-r8}
\end{figure}

Figure~\ref{fig:final-chair-r8} also shows a qualitative difference between the
sigmoid models and the direct tensor baselines. At the same nominal rank, the
sigmoid reconstructions retain small and thin occupied components of the chair,
whereas the direct baselines tend to smooth them out or omit them. The random
CE-Tucker reconstruction is visibly degraded.

This behavior is consistent with
the initialization study: in the stored Chair and Dragon runs, random
CE-Tucker was systematically worse than the structured starts at the same
rank, including the repeated $512^3$ runs. It should therefore be interpreted
as an initialization--optimization effect, rather than as evidence that the
Tucker parameterization has insufficient representational capacity. 

The visual comparison illustrates why overlap metrics are useful in addition to
cross-entropy. Methods with larger IoU preserve the occupied surface more
continuously, whereas lower-IoU methods lose thin components or introduce
disconnected fragments. The Dragon reconstructions at $r=4, 8$ and Chair at $r=4$ rare provided in Appendix~\ref{sec:appendix-visual}.



\section{Arithmetic complexity}
\label{sec:complexity}

We compare the computational cost of the methods used in the experiments by
counting multiplication-dominated arithmetic operations on a dense
third-order tensor of size $n\times n\times n$. All Tucker and
TT ranks are denoted by $r$, and $N_{\mathrm{iter}}$ denotes the number of
optimization iterations. 

The total costs are summarized in
Table~\ref{tab:algorithm-complexity}. We denote the selected one-time
initialization costs by $I_{\mathrm{init}}^{\mathrm{Tucker}}$ and
$I_{\mathrm{init}}^{\mathrm{TT}}$, which depend on the chosen strategy. The
alternatives are listed in
Table~\ref{tab:initialization-complexity}.

\begin{table}[h]
\centering
\caption{Arithmetic complexity of the compared algorithms for an
$n\times n\times n$ tensor with rank parameter $r$.}
\label{tab:algorithm-complexity}
\begin{tabular}{p{0.23\linewidth}p{0.29\linewidth}p{0.38\linewidth}}
\hline
\textbf{Algorithm} & \textbf{Initialization cost} & \textbf{Algorithm cost}\\
\hline
TT-TSDF & $-$ & $\mathcal{O}(n^3r+n^2r^3)$ (TT-SVD)\\
Tucker-binary & $-$ & $\mathcal{O}(n^3r)$ (ST-HOSVD)\\
TT-binary & $-$ & $\mathcal{O}(n^3r+n^2r^3)$ (TT-SVD)\\
CE-Tucker & $I_{\mathrm{init}}^{\mathrm{Tucker}}$ &
$\mathcal{O}\!\left(N_{\mathrm{iter}}(n^3r+nr^2)\right)$\\
CE-TT & $I_{\mathrm{init}}^{\mathrm{TT}}$ &
$\mathcal{O}\!\left(N_{\mathrm{iter}}(n^3r^2+nr^3)\right)$\\
\hline
\end{tabular}
\end{table}
\vspace{-1cm}
For the one-shot baselines, ST-HOSVD
computes a truncated SVD of each mode unfolding and projects the tensor onto
the resulting subspaces, yielding the dense third-order estimate
$\mathcal{O}(n^3r)$ in our accounting. TT-SVD successively reshapes the tensor,
computes a truncated SVD, and passes the resulting singular values to the next
core with final complexity as 
$\mathcal{O}(n^3r+n^2r^3)$. These are standard dense costs for ST-HOSVD and
TT-SVD~\cite{delathauwer2000,oseledets2011tt}. For CE-Tucker, one block-gradient update contracts the dense error tensor with
the fixed factor matrices and core, giving $\mathcal{O}(n^3r)$ operations with
QR orthogonalization of a factor costs $\mathcal{O}(nr^2)$. For CE-TT, the
left and right environments of an active core give the corresponding costs
$\mathcal{O}(n^3r^2)$ for the gradient and $\mathcal{O}(nr^3)$ for TT
canonicalization. Multiplication by $N_{\mathrm{iter}}$ gives the totals in
Table~\ref{tab:algorithm-complexity}. 

For initialization, let $r_0<r$ be the lower rank used by a rank-expansion or
residual-informed start and define $q=r-r_0$. We introduce the following
notation to simplify the initialization-complexity expressions:
\begin{align*}
I_{\mathrm{rand}}^{\mathrm{Tucker}}(r)&=\mathcal{O}(nr+r^3),
&I_{\mathrm{rand}}^{\mathrm{TT}}(r)&=\mathcal{O}(nr^2),\\
I_{\mathrm{dec}}^{\mathrm{Tucker}}(r)&=\mathcal{O}(n^3r),
&I_{\mathrm{dec}}^{\mathrm{TT}}(r)&=\mathcal{O}(n^3r+n^2r^3),\\
I_{\mathrm{res}}^{\mathrm{Tucker}}(q)&=\mathcal{O}(n^3r_0+n^3q),
&I_{\mathrm{res}}^{\mathrm{TT}}(q)&=\mathcal{O}(n^3r_0^2+n^3q+n^2q^3).
\end{align*}
\vspace{-1cm}
\begin{table}[h]
\centering
\caption{Arithmetic costs of the four initialization strategies.}
\label{tab:initialization-complexity}
\scriptsize
\setlength{\tabcolsep}{0pt}
\begin{tabular}{@{}p{0.09\linewidth}p{0.14\linewidth}p{0.20\linewidth}p{0.29\linewidth}p{0.29\linewidth}@{}}
\hline
 & \textbf{Random} & \mbox{\textbf{Decomposition}} & \mbox{\textbf{Rank-expansion}} &
\mbox{\textbf{Residual}}\\
\hline
CE-Tucker & $I_{\mathrm{rand}}^{\mathrm{Tucker}}(r)$ & $I_{\mathrm{dec}}^{\mathrm{Tucker}}(r)$ &
$I_{\mathrm{dec}}^{\mathrm{Tucker}}(r_0)+I_{\mathrm{rand}}^{\mathrm{Tucker}}(q)$ &
$I_{\mathrm{dec}}^{\mathrm{Tucker}}(r_0)+I_{\mathrm{res}}^{\mathrm{Tucker}}(q)$\\
CE-TT & $I_{\mathrm{rand}}^{\mathrm{TT}}(r)$ & $I_{\mathrm{dec}}^{\mathrm{TT}}(r)$ &
$I_{\mathrm{dec}}^{\mathrm{TT}}(r_0)+I_{\mathrm{rand}}^{\mathrm{TT}}(q)$ &
$I_{\mathrm{dec}}^{\mathrm{TT}}(r_0)+I_{\mathrm{res}}^{\mathrm{TT}}(q)$\\
\hline
\end{tabular}
\normalsize
\end{table}

A random start generates and stores the factors and cores. For the
third-order tensor considered here, this costs $I_{\mathrm{rand}}^{\mathrm{Tucker}}(r)$ for
Tucker and $I_{\mathrm{rand}}^{\mathrm{TT}}(r)$ for TT. For a general
$d$-way tensor, the corresponding costs are $\mathcal{O}(dnr+r^d)$ and
$\mathcal{O}(dnr^2)$, respectively. A decomposition start applies ST-HOSVD
for Tucker or TT-SVD for TT, whose costs are
$I_{\mathrm{dec}}^{\mathrm{Tucker}}(r)$ and
$I_{\mathrm{dec}}^{\mathrm{TT}}(r)$, respectively. For a rank-expansion warm start, the newly added
random components have the same costs as a random start with rank $q$, while
the lower-rank decomposition contributes the corresponding
$I_{\mathrm{dec}}^{\mathrm{Tucker}}(r_0)$ or
$I_{\mathrm{dec}}^{\mathrm{TT}}(r_0)$. No QR
orthogonalization or TT canonicalization is included in these estimates. A
residual-informed start additionally reconstructs the lower-rank approximation
and decomposes the residual, as represented by the corresponding
$I_{\mathrm{res}}^{\mathrm{Tucker}}(q)$ or
$I_{\mathrm{res}}^{\mathrm{TT}}(q)$.
Dense tensor formation and residual subtraction are excluded, whereas generation
of the newly added random components is included. For the
decomposition start, exact SVD can be
replaced by randomized range-SVD with a small oversampling parameter and a
fixed number of power iterations. This avoids forming a full dense SVD and
reduces the practical initialization cost, while preserving the same leading
rank-truncated order shown in the table.

\section{Conclusion}\label{sec:conclusion}

We introduced nonlinear tensor decomposition for binary data using sigmoid Tucker and tensor-train representations minimized with binary cross-entropy. The complete optimization problem is jointly non-convex, but it is convex with respect to one factor or core block when the remaining blocks are fixed. The experiments show that both rank and initialization influence the final reconstruction. In particular, at small ranks an informative initialization can improve IoU and Hamming error relative to the standard direct tensor baselines, while increasing the rank improves the attainable reconstruction quality.

Several directions remain open. First, a more systematic convergence analysis is
needed to explain the dependence of block-coordinate optimization on the tensor
format, rank, and initialization. This analysis should also lead to principled
procedures for selecting the rank and initialization from the data, rather than
relying on a fixed rank grid and manually specified starts. Rank-adaptive
stopping rules and quantized factor storage are further important extensions.
Second, it is worthwhile to characterize nonlinear functions that admit a
separable representation, or an analogous tensor-valued factorization, so that
their evaluation can be implemented with efficient matrix--vector products,
elementwise operations, and tensor contractions. Such representations could
reduce the cost of nonlinear tensor multiplication while retaining the observed
reconstruction advantages. Finally, the proposed framework should be extended
to other data formats, including bounded real-valued, count, categorical, and
multimodal observations, with the nonlinear map and discrepancy chosen to match
the corresponding observation model.

\appendix
\section{Additional visual reconstructions}\label{sec:appendix-visual}

\begin{figure}[!h]
\centering
\includegraphics[width=0.85\textwidth]{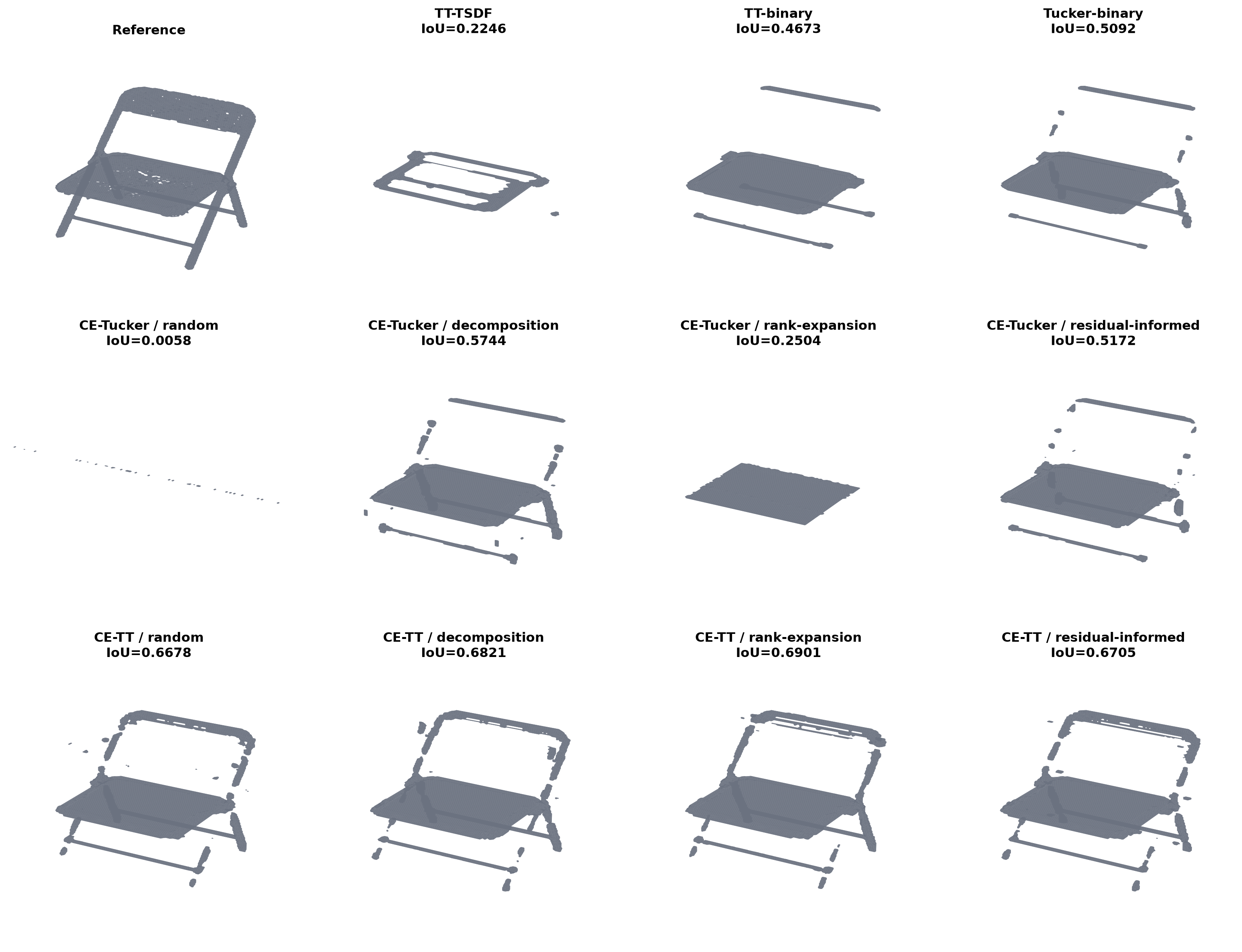}
\caption{ModelNet40 Chair at $512^3$ and $r=4$.}
\label{fig:appendix-chair-r4}
\end{figure}

\begin{figure}[!h]
\centering
\includegraphics[width=0.85\textwidth]{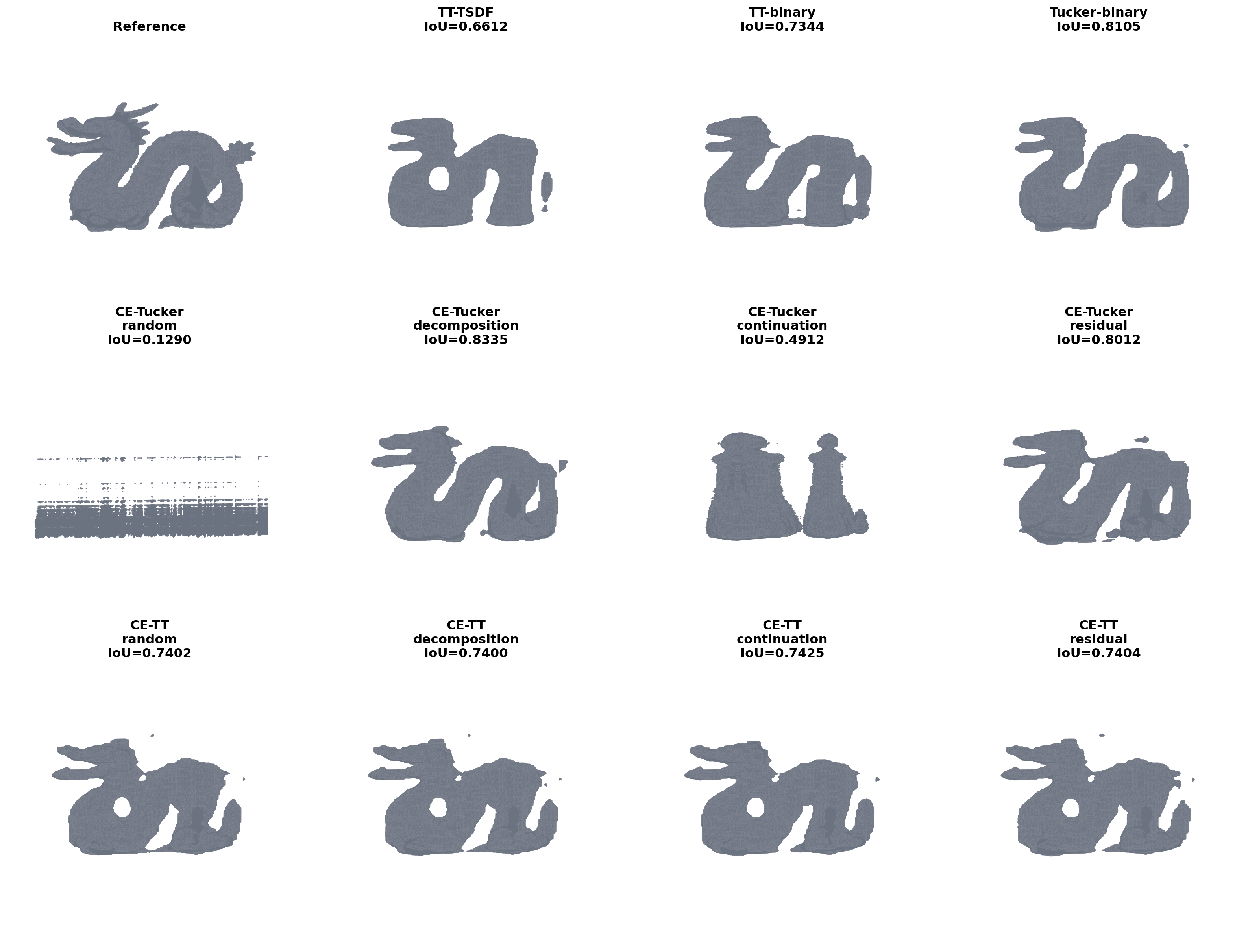}
\caption{Stanford Dragon at $512^3$ and $r=4$.}
\label{fig:appendix-dragon-r4}
\end{figure}

\begin{figure}[!h]
\centering
\includegraphics[width=0.85\textwidth]{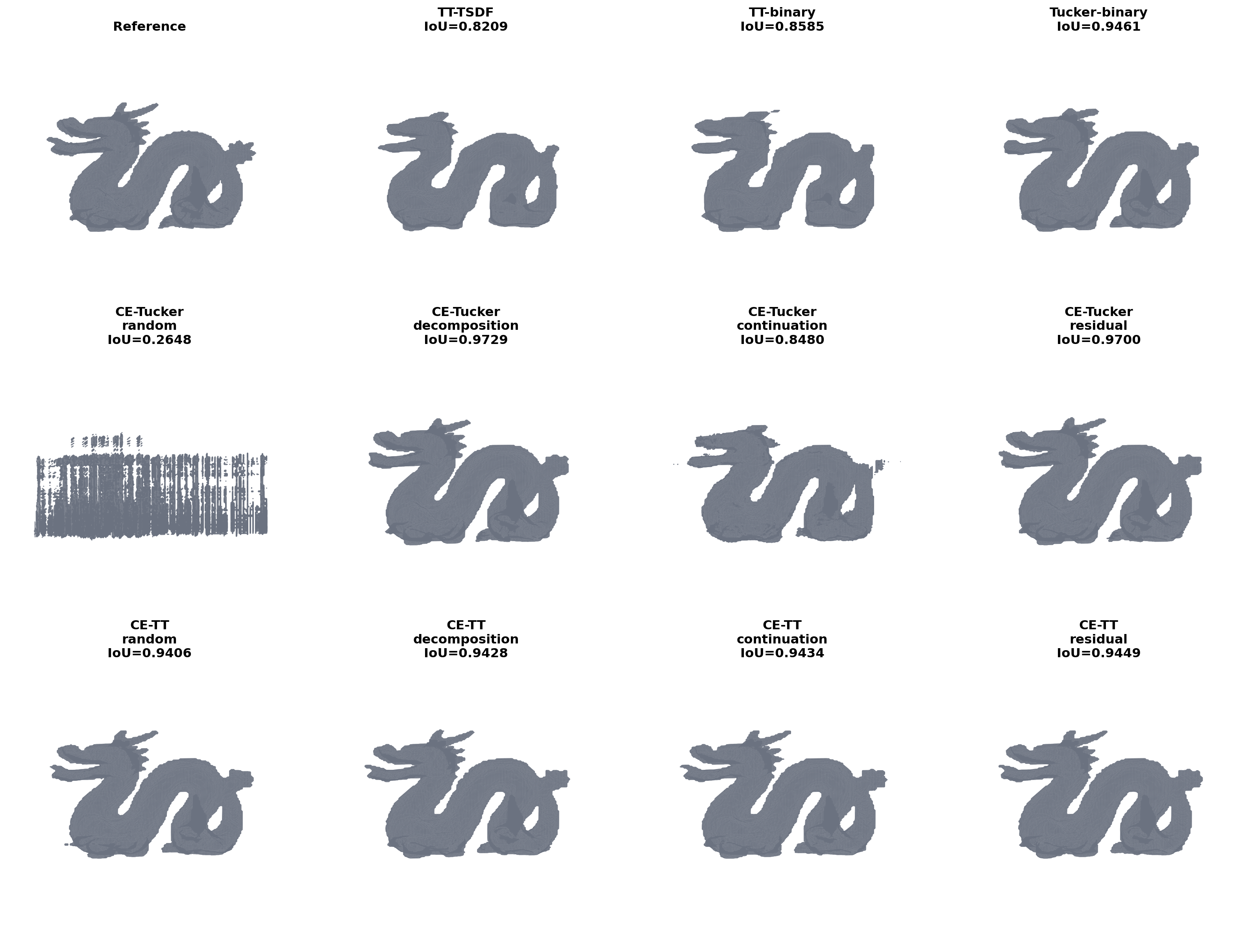}
\caption{Stanford Dragon at $512^3$ and $r=8$.}
\label{fig:appendix-dragon-r8}
\end{figure}

For additional visual comparisons, Fig.~\ref{fig:appendix-chair-r4} shows the
Chair at $r=4$, while Figs.~\ref{fig:appendix-dragon-r4} and
\ref{fig:appendix-dragon-r8} show the Stanford Dragon at $r=4$ and $r=8$,
respectively. 

\newpage


\bibliography{refs}

@article{eckart1936,
  author = {Eckart, Carl and Young, Gale},
  title = {The Approximation of One Matrix by Another of Lower Rank},
  journal = {Psychometrika}, year = {1936}, volume = {1}, number = {3}, pages = {211--218},
  doi = {10.1007/BF02288367}
}

@article{von1937some,
  title={Some matrix-inequalities and metrization of matric-space},
  author={von Veumann, J},
  journal={Tomsk Univ Rev},
  volume={1},
  pages={286--300},
  year={1937}
}

@article{mirsky1960,
  author = {Mirsky, Leon},
  title = {Symmetric Gauge Functions and Unitarily Invariant Norms},
  journal = {The Quarterly Journal of Mathematics},
  year = {1960},
  volume = {11},
  number = {1},
  pages = {50--59},
  doi = {10.1093/qmath/11.1.50}
}

@article{lee1999nmf,
  author = {Lee, Daniel D. and Seung, H. Sebastian},
  title = {Learning the Parts of Objects by Non-Negative Matrix Factorization},
  journal = {Nature},
  year = {1999},
  volume = {401},
  pages = {788--791},
  doi = {10.1038/44565}
}

@inproceedings{chierichetti2017l1,
  author = {Chierichetti, Flavio and Gollapudi, Sreenivas and Kumar, Ravi and Lattanzi, Silvio and Panigrahy, Rina and Woodruff, David P.},
  title = {Algorithms for $\ell_p$ Low-Rank Approximation},
  booktitle = {Proceedings of the 34th International Conference on Machine Learning},
  year = {2017},
  volume = {70},
  pages = {806--814},
  month = {06--11 Aug},
  editor = {Precup, Doina and Teh, Yee Whye},
  publisher = {PMLR},
  series = {Proceedings of Machine Learning Research},
  address = {Sydney, Australia},
  url = {https://proceedings.mlr.press/v70/chierichetti17a.html}
}

@inproceedings{lee2001nmf,
  author = {Lee, Daniel D. and Seung, H. Sebastian},
  title = {Algorithms for Non-Negative Matrix Factorization},
  booktitle = {Advances in Neural Information Processing Systems},
  year = {2001},
  volume = {13},
  pages = {556--562},
  publisher = {MIT Press},
  address = {Cambridge, MA, USA},
  url = {https://proceedings.neurips.cc/paper/1861-algorithms-for-non-negative-matrix-factorization}
}

@article{hoyer2004sparse,
  author = {Hoyer, Patrik O.},
  title = {Non-Negative Matrix Factorization with Sparseness Constraints},
  journal = {Journal of Machine Learning Research},
  year = {2004},
  volume = {5},
  pages = {1457--1469},
  url = {https://www.jmlr.org/papers/v5/hoyer04a.html}
}

@article{mairal2010online,
  author = {Mairal, Julien and Bach, Francis and Ponce, Jean and Sapiro, Guillermo},
  title = {Online Learning for Matrix Factorization and Sparse Coding},
  journal = {Journal of Machine Learning Research},
  year = {2010},
  volume = {11},
  pages = {19--60},
  url = {https://www.jmlr.org/papers/v11/mairal10a.html}
}

@book{bishop2006pattern,
  author = {Bishop, Christopher M.},
  title = {Pattern Recognition and Machine Learning},
  publisher = {Springer},
  address = {New York, NY, USA},
  year = {2006},
  doi = {10.1007/978-0-387-45528-0},
  url = {https://doi.org/10.1007/978-0-387-45528-0}
}

@article{delathauwer2000,
  author = {De Lathauwer, Lieven and De Moor, Bart and Vandewalle, Joos},
  title = {A Multilinear Singular Value Decomposition},
  journal = {SIAM Journal on Matrix Analysis and Applications},
  year = {2000}, volume = {21}, number = {4}, pages = {1253--1278},
  doi = {10.1137/S0895479896305696}
}

@article{kolda2009,
  author = {Kolda, Tamara G. and Bader, Brett W.},
  title = {Tensor Decompositions and Applications},
  journal = {SIAM Review}, year = {2009}, volume = {51}, number = {3}, pages = {455--500},
  doi = {10.1137/07070111X}
}

@article{oseledets2011tt,
  author = {Oseledets, Ivan V.},
  title = {Tensor-Train Decomposition},
  journal = {SIAM Journal on Scientific Computing}, year = {2011}, volume = {33}, number = {5},
  pages = {2295--2317}, doi = {10.1137/090752286}
}

@article{holtz2012als,
  author = {Holtz, Sebastian and Rohwedder, Thorsten and Schneider, Reinhold},
  title = {The Alternating Linear Scheme for Tensor Optimization in the Tensor Train Format},
  journal = {SIAM Journal on Scientific Computing},
  year = {2012}, volume = {34}, number = {2}, pages = {A683--A713},
  doi = {10.1137/100818893}
}

@article{rohwedder2013als,
  author = {Rohwedder, Thorsten and Uschmajew, Andr{\'e}},
  title = {On Local Convergence of Alternating Schemes for Optimization of Convex Problems in the Tensor Train Format},
  journal = {SIAM Journal on Numerical Analysis},
  year = {2013}, volume = {51}, number = {2}, pages = {1134--1162},
  doi = {10.1137/110857520}
}

@article{chewei2019randomized,
  author = {Che, Maolin and Wei, Yimin},
  title = {Randomized Algorithms for the Approximations of Tucker and the Tensor Train Decompositions},
  journal = {Advances in Computational Mathematics},
  year = {2019},
  volume = {45},
  number = {1},
  pages = {395--428},
  doi = {10.1007/s10444-018-9622-8}
}

@article{saul2022,
  author = {Saul, Lawrence K.},
  title = {A Nonlinear Matrix Decomposition for Mining the Zeros of Sparse Data},
  journal = {SIAM Journal on Mathematics of Data Science}, year = {2022}, volume = {4}, number = {2},
  pages = {431--463}, doi = {10.1137/21M1405769},
}

@inproceedings{seraghiti2023,
  author = {Seraghiti, Giovanni and Awari, Atharva and Vandaele, Arnaud and Porcelli, Margherita and Gillis, Nicolas},
  title = {Accelerated Algorithms for Nonlinear Matrix Decomposition with the ReLU Function},
  booktitle = {2023 IEEE 33rd International Workshop on Machine Learning for Signal Processing (MLSP)},
  year = {2023}, pages = {1--6}, publisher = {IEEE},
  address = {Rome, Italy},
  doi = {10.1109/MLSP55844.2023.10285984}
}

@article{wang2024,
  author = {Wang, Qingsong and Cui, Chunfeng and Han, Deren},
  title = {A Momentum Accelerated Algorithm for ReLU-Based Nonlinear Matrix Decomposition},
  journal = {IEEE Signal Processing Letters}, year = {2024}, volume = {31}, pages = {2865--2869},
  doi = {10.1109/LSP.2024.3475910}
}

@article{gillis2025,
  author = {Gillis, Nicolas and Porcelli, Margherita and Seraghiti, Giovanni},
  title = {An Extrapolated and Provably Convergent Algorithm for Nonlinear Matrix Decomposition with the ReLU Function},
  journal = {SIAM Journal on Mathematics of Data Science}, year = {2026}, volume = {8}, number = {3},
  pages = {738--765},
  doi = {10.1137/25M1746859},
  url = {https://doi.org/10.1137/25M1746859}
}

@article{wang2025,
  author = {Wang, Qingsong and Qu, Yunfei and Cui, Chunfeng and Han, Deren},
  title = {An Accelerated Alternating Partial Bregman Algorithm for ReLU-Based Matrix Decomposition},
  journal = {Journal of Scientific Computing}, year = {2025}, volume = {105}, number = {2}, pages = {48},
  doi = {10.1007/s10915-025-03067-w}
}

@inproceedings{nguyen2025,
  author = {Nguyen, Harrison and Awari, Atharva Abhijit and Vandaele, Arnaud and Gillis, Nicolas},
  title = {Nonlinear Matrix Decomposition with the Sigmoid Function},
  booktitle = {2025 IEEE 35th International Workshop on Machine Learning for Signal Processing (MLSP)},
  year = {2025},
  publisher = {IEEE},
  address = {Istanbul, Turkey},
  doi = {10.1109/MLSP62443.2025.11204338}
}

@inproceedings{thayse2026,
  author = {Thayse, Alex and Vandaele, Arnaud and Gillis, Nicolas},
  title = {Sigmoid Nonlinear Matrix Decomposition with Cross Entropy},
  booktitle = {2026 34th European Signal Processing Conference (EUSIPCO)},
  year = {2026},
  pages = {2551--2551},
  address = {Bruges, Belgium},
  publisher = {IEEE},
  isbn = {978-94-645936-3-1},
  url = {https://eurasip.org/Proceedings/Eusipco/Eusipco2026/HTML/author-index/index.html}
}

@misc{awari2025admm,
  author = {Awari, Atharva and Gillis, Nicolas and Vandaele, Arnaud},
  title = {Alternating Direction Method of Multipliers for Nonlinear Matrix Decompositions},
  year = {2025},
  eprint = {2512.17473},
  archivePrefix = {arXiv},
  primaryClass = {eess.SP},
  url = {https://arxiv.org/abs/2512.17473}
}

@misc{awari2026ncp,
  author = {Awari, Atharva and van der Ploeg, Geert Roelof and Vandaele, Arnaud and Gillis, Nicolas and Acar, Evrim},
  title = {Nonlinear Tensor Decomposition for Pattern Discovery},
  year = {2026},
  eprint = {2610.05339},
  archivePrefix = {arXiv},
  primaryClass = {cs.LG},
  doi = {10.48550/arXiv.2610.05339},
  url = {https://arxiv.org/abs/2610.05339}
}

@article{sokolova2009metrics,
  author = {Sokolova, Marina and Lapalme, Guy},
  title = {A Systematic Analysis of Performance Measures for Classification Tasks},
  journal = {Information Processing \& Management},
  year = {2009},
  volume = {45},
  number = {4},
  pages = {427--437},
  doi = {10.1016/j.ipm.2009.03.002}
}

@article{wang2026relu_ntd,
  author = {Wang, Qingsong and Liu, Zehui and Cui, Chunfeng and Han, Deren},
  title = {ReLU-Based Nonlinear Tensor Decomposition for Sparse Nonnegative Data},
  journal = {Journal of Computational and Applied Mathematics}, year = {2027}, volume = {491},
  pages = {118016},
  doi = {10.1016/j.cam.2026.118016},
  url = {https://doi.org/10.1016/j.cam.2026.118016},
  note = {Published online in 2026; volume 491 (2027)}
}

@inproceedings{boyko2020,
  author = {Boyko, Alexey I. and Matrosov, Mikhail P. and Oseledets, Ivan V. and Tsetserukou, Dzmitry and Ferrer, Gonzalo},
  title = {{TT-TSDF}: Memory-Efficient {TSDF} with Low-Rank Tensor Train Decomposition},
  booktitle = {2020 IEEE/RSJ International Conference on Intelligent Robots and Systems (IROS)},
  year = {2020},
  pages = {5706--5711},
  doi = {10.1109/IROS45743.2020.9341464}
}

@article{DolgovSavostyanov2014,
  author  = {Dolgov, Sergey V. and Savostyanov, Dmitry V.},
  title   = {Alternating Minimal Energy Methods for Linear Systems in Higher Dimensions},
  journal = {SIAM Journal on Scientific Computing},
  volume  = {36},
  number  = {5},
  pages   = {A2248--A2271},
  year    = {2014},
  doi     = {10.1137/140953289}
}

\end{document}